\documentclass[12pt,oneside]{amsart}

\usepackage{amsfonts,amssymb}

\usepackage{a4wide}
\usepackage{manfnt,mathabx}

\usepackage{comment}

\def\BS{\mathrm{BS}}
\def\pf{{\noindent \it Proof. \ }}
\def\qed{\hfill $\Box$}
\def\g{\gamma}
\def\f{\phi}
\def\ar{\rightarrow}
\def\<{\langle}
\def\>{\rangle}
\def\N{\mathbb N}
\def\Z{\mathbb Z}

\def\f{\phi}
\def\e{\varepsilon}
\def\w{\omega}
\newtheorem{lemma}{Lemma}
\newtheorem{theorem}{Theorem}
\newtheorem{cor}{Corollary}
\newtheorem{proposition}{Proposition}
\newtheorem{remark}{Remark}

\reversemarginpar

\title{Baumslag-Solitar groups and residual nilpotence} 
\author{C.E.  Kofinas, V.  Metaftsis and A.I.  Papistas}
\address{Department of Mathematics, University of the Aegean,  Karlovassi,  832 00  Samos, Greece}\email{kkofinas@aegean.gr}
\address{Department of Mathematics, University of the Aegean,  Karlovassi, 832 00  Samos, Greece}\email{vmet@aegean.gr}
\address{Department of Mathematics,  Aristotle University of Thessaloniki,  541 24 Thessaloniki, Greece.}\email{apapist@math.auth.gr}
\date{}

\begin{document}

\begin{abstract}
Let $G$ be a Baumslag-Solitar group. We calculate the intersection $\g_{\w}(G)$ of all terms of the lower central series of $G$.  Using this, we show that $[\g_{\w}(G),G]=\g_{\w}(G)$; thus answering a question of Bardakov and Neschadim \cite{bn}.  For any $c \in \mathbb{N}$, with $c \geq 2$, we show,  by using Lie algebra methods, that the quotient group $\gamma_{c}(G)/\gamma_{c+1}(G)$ of the lower central series of $G$ is finite.
\end{abstract}

\maketitle
\section{Introduction}

Baumslag-Solitar groups are groups that admit a presentation of the form 
$$\BS(m,n)=\langle a, t \mid t^{-1}a^{m}t = a^{n}\rangle$$ where $m, n$ are non-zero integers. They were  introduced in \cite{baso} as examples of two-generator one-relator groups with proper quotients isomorphic to the group itself (that is, the groups don't satisfy the Hopf property).  Since then,  Baumslag-Solitar groups and their properties have been extensively studied by various authors
and they have been the testbed for various conjectures and theories.  

Our work is mainly concerned with the residual nilpotence of these groups.   A survey about the residual properties of these groups is given in \cite{mold-rev}.   In \cite{bn}, Bardakov and Neschadim studied the lower central series of Baumslag-Solitar groups and computed the intersection of all terms of the lower central series for some special cases of the non-residually nilpotent Baumslag-Solitar groups. In the present paper,  one of our aims is to compute explicitly the intersection of all terms of the lower central series for any Baumslag-Solitar group which is not residually nilpotent.

 Throughout this paper, a Baumslag-Solitar group is denoted by $\BS(m, n)$. Since $\BS(m, n)$, $\BS(n, m)$ and $\BS(-m, -n)$ are pairwise isomorphic, we may assume, without loss of generality, that the integers $m$ and $n$ in the presentation of $\BS(m, n)$ satisfy the condition $0 < m \leq |n|$. 

One of our main results is the following.

\begin{theorem}\label{main}
Let $G=\BS(m,n)$, with $0 < m \leq |n|$, ${\rm gcd}(m,n)=d \geq 1$, $m=dm_1$ and $n=dn_1$.  
\begin{enumerate}

\item If $n_1\not\equiv m_1 \pmod p$ for every prime number $p$, then $\g_{\w}(G)$ is the normal closure of the set $\{a^{d}, [t^{-k}a^{\mu}t^{k},a^{\nu}]\mid k\in\Z, \mu,\nu \in \mathbb{N}, {\rm gcd}(\mu, \nu) = 1 ~{\rm and}~\mu \nu = d\}$ in $G$. 

\item If there is a prime number $p$ such that  $n_1\equiv m_1\pmod p$, then $\g_{\w}(G)$ is the normal closure  of the set $\{[t^{-k}a^{\mu}t^{k},a^{\nu}]\mid k\in\Z, \mu,\nu \in \mathbb{N}, {\rm gcd}(\mu, \nu) = 1 ~{\rm and}~\mu \nu = d\}$ in $G$.      

\end{enumerate}
\end{theorem}

Next,  we are concerned with the following question of \cite{bn}.  Let $G={\rm BS}(m,n)$, with $0< m\le |n|$.  Is it true that $\gamma_{\w}(G) = [\gamma_{\w}(G), G]$?
In fact, we are able to answer the above question affirmatively.  

\begin{theorem}\label{thefin}
Let $G={\rm BS}(m,n)$, with $0 < m \leq |n|$. Then $[\g_{\w}(G),G]=\g_{\w}(G)$.     
\end{theorem}

In section $5$,  by using Lie algebra methods, we show that for a Baumslag-Solitar group $G$ the quotient groups $\g_c(G)/\g_{c+1}(G)$, with $c \geq 2$, of the lower central series of $G$ are finite.

\section{Auxiliary results}

Let $G$ be a group.  For elements $a, b$ of $G$, we write $[a, b]$ for the commutator of $a$ and $b$, that is $[a, b] = a^{-1}b^{-1}ab$. 
We denote $\langle g_{1}, \ldots, g_{c} \rangle$ the subgroup of $G$ generated by the elements $g_{1}, \ldots, g_{c}$. For subgroups $A$ and $B$ of $G$, we write $[A, B] = \langle [a, b], a \in A, b \in B \rangle$. For a positive integer $c$, let $\g_1(G)=G$ and, for $c\ge 2$, let $\g_c(G)=[\g_{c-1}(G),G]$ be the $c$-th term of the lower central series of $G$. We point out that $\gamma_{2}(G) = [G,G] = G^{\prime}$; that is, the derived group of $G$. We write $\gamma_{\w}(G)$ for the intersection of all terms of the lower central series of $G$, that is, $\gamma_{\w}(G) = \bigcap_{c\geq 1}\gamma_{c}(G)$. We say $G$ is a residually $\mathcal{P}$ group if for every element $1 \neq g \in G$ there is a normal subgroup $N_{g}$ of $G$ not containing $g$ such that $G/N_{g}$ has the property $\mathcal{P}$.  In case ${\mathcal P}$ is nilpotency, we say that the group is residually nilpotent.  Equivalently,   we say that $G$ is a residually nilpotent group, if $\gamma_{\w}(G) = \{1\}$. For the rest of the paper,  $N_{\w}$ denotes the intersection of all finite index normal  subgroups of $G$ and $(Np)_{w}$ denotes the intersection of all finite index normal subgroups of $G$ with index some power of a prime number $p$.

The following proposition summarizes some residual properties concerning Baumslag-Solitar groups.

\begin{proposition}\label{known}
Let $G$ be the Baumslag-Solitar group with presentation
$$G=\BS(m,n)=\<t,a \mid t^{-1}a^mt=a^n\>,$$
with $0 < m \leq |n|$. Then,
\begin{enumerate}
\item The group $G$ is residually finite if and only if $m=1$ or $|n|=m$.

\item The group $G$ is residually nilpotent if and only if $m=1$ and $n\neq 2$ or $|n|=m=p^r$, $r > 0$ for some prime number $p$.

\item The group $G$ is residually finite $p$-group for some prime number $p$ if and only if $m=1$ and $n\equiv 1\pmod p$ or $n=m$ and $m=p^r$ or $n=-m$, $p=2$ and $m=2^r$, $r\ge 1$.
\end{enumerate}
\end{proposition}

\begin{remark}\upshape
\begin{enumerate}
\item The residual finiteness of the Baumslag-Solitar groups originally studied in \cite{baso} and completed in \cite{meskin}. Recently,  Moldavanskii in \cite{mold} calculated $N_{\w}$ for Baumslag-Solitar groups.

\item In \cite{rava}, Raptis and Varsos gave necessary conditions for the residual
nilpotence of HNN-extensions with base group a finitely generated abelian group. Proposition \ref{known} (1) follows from \cite{rava}: it is a special case of \cite[Corollary 2.7]{rava}. 

\item Proposition \ref{known} (3) follows from the study of Kim and McCarron in certain one relator groups (see \cite[Main Theorem]{kimc}).  Also, Moldavanskii in \cite{mold-ru} (see also \cite{mold-rev})  calculated $(N_p)_{\w}$ for $\BS(m,n)$.
\end{enumerate}
\end{remark}

\begin{lemma}\label{rep5}
Let $G$ be a group and $N$ be a normal subgroup of $G$.  Let $x, y \in G$ such that $[x, y] \in N$.
\begin{enumerate}
\item Then for all $\kappa\in\N$, $[x^{\kappa}, y],  [x,y^{\kappa}] \in N$ and $[x^{\kappa}, y] \equiv [x,y^{\kappa}] \equiv [x,y]^{\kappa} ~({\rm mod}~[N,G])$.

\item If $[x, y^{m}] \in [N, G]$ for some $m \in \mathbb{N}$, then $[x, y]^{m} \in [N, G]$.

\item If $[x^{m}, y] \in [N, G]$ for some $m \in \mathbb{N}$, then $[x, y]^{m} \in [N, G]$.
\end{enumerate}
\end{lemma}

\pf 
\begin{enumerate}

\item This is straightforward. 

\item Let $[x,y^{m}] \in [N,G]$ for some $m \in \mathbb{N}$. By Lemma \ref{rep5}~(1) (for $\kappa = m$), we have $[x,y^{m}] = [x,y]^{m}w$, with $w \in [N,G]$, and so, $[x,y]^{m} \in [N,G]$.   

\item Let $[x^{m},y] \in [N,G]$. By Lemma \ref{rep5}~(1) (for $\kappa = m$), we have $[x^{m},y] = [x,y]^{m}w_{1}$, with $w_{1} \in [N,G]$, and so, $[x,y]^{m} \in [N,G]$.  \qed
\end{enumerate}

The following result gives us a relation among residually finite, residually nilpotent and residually finite $p$-group for some prime number $p$.

\begin{lemma}\label{bounds}
Let $G$ be a finitely generated group. Then $N_{\w}\le\g_{\w}(G)\le (Np)_{w}$. Moreover, $\bigcap\limits_{p\ {\rm prime}}(Np)_{w}\le \g_{\w}(G)$.
\end{lemma}

\pf Since a finite $p$-group is nilpotent, we have every residually finite $p$-group is also residually nilpotent. Hence $G/(Np)_{\omega}$ is residually nilpotent and so $\gamma_{\omega}(G)\le (Np)_{\omega}$.  Since $G$ is finitely generated and $G/\gamma_{\omega}(G)$ is residually nilpotent, we have $G/\gamma_{\omega}(G)$ is  residually finite.  We claim that $N_{\omega}\le \gamma_{\omega}(G)$.  Let $g \in N_{\omega}$ and $g \notin \gamma_{\omega}(G)$. Since $G/\gamma_{\omega}(G)$ is  residually finite,  there exists a normal subgroup $N_{g}$ of $G$ such that $g \notin N_{g}$, $\gamma_{\omega}(G) \subseteq N_{g}$ and $G/N_{g}$ is finite, which is a contradiction since $g \in N_{\omega}$.

Write $B = \bigcap\limits_{p\ {\rm prime}}(Np)_{\omega}$. To get a contradiction, we assume that $g \in B$ and $g \notin \gamma_{\omega}(G)$.  In the next few lines, let $\widetilde{G}=G/\gamma_{\omega}(G)$. Since $G$ is finitely generated and $\widetilde{G}$ is residually nilpotent, there exists an epimorphism $\phi$ from $\widetilde{G}$ onto a finitely generated nilpotent group $H$ with $\phi(g\gamma_{\omega}(G)) \neq 1$. Since $H$ is polycyclic, we have $H$ is residually finite.  Thus there exists a finite nilpotent group $\hat{H}$ and an epimorphism $\hat{\phi}: \widetilde{G} \rightarrow \hat{H}$ with $\hat{\phi}(g \gamma_{\omega}(G)) \neq 1$. Since $\hat{H}$ is the direct product of its Sylow $p$-subgroups, there exist a prime number $p$ and a Sylow $p$-subgroup $S_{p}$ of $\hat{H}$ such that $\hat{\phi}(g \gamma_{\omega}(G)) \in S_{p} \setminus \{1\}$. Since $S_{p}$ is a finite $p$-group, we have $g \notin (Np)_{\omega}$, and so $g\not\in B$,  which is a contradiction. \qed

\begin{lemma}\label{zzm}
Let 
$$G=\<t,x_1,\ldots, x_n\mid x_i^{p_i^{r_i}}=1,  i=1,\ldots, n, [x_i,x_j]=1 \>\cong \Z* \left(\Z_{p_1^{r_1}}\times\dots\times\Z_{p_n^{r_n}}\right)$$ where $n\ge 2$, $p_{1}, \ldots, p_{n}$ are distinct prime numbers. Then $\g_{\w}(G)$ is the normal closure of the set  $\{[t^{-k} x_it^{k}, x_j]: i,j\in\{1,\ldots, n\}; i\neq j; k\in\Z\}$ in $G$. 
\end{lemma}

\pf  The elements $x_i$ have orders $p_i^{r_i}$ \and so  the orders of $x_i$ and $x_j$ are coprime for every $i\neq j$.  
We write $\widetilde{G} = G/\gamma_{\omega}(G)$.  As in the proof of Lemma \ref{bounds}, there exist a finite nilpotent group $\hat{H}$ and an epimorphism $\hat{\phi}: \widetilde{G}\rightarrow\hat{H}$ with $\hat{\phi}({\widetilde{g}})\neq 1$ for every $1\neq\widetilde{g}\in\widetilde{G}$. 
Furthermore there exists a prime number $p$ and a Sylow $p$-subgroup $S_p$ of $\hat{H}$ such that $\hat{\phi}({\widetilde{g}})\in S_p\setminus\{1\}$.  Hence the elements $[t^{-k}x_{i}t^{k},x_{j}]$, with $i \neq j$, are trivial under any homomorphism onto $S_{p}$.  Therefore $[t^{-k}x_{i}t^{k},x_{j}] \in \gamma_{\omega}(G)$ for all $k \in \mathbb{Z}$ and $i, j \in \{1, \ldots, n\}$, with $i \neq j$. 

For the converse, we will show that $G/N$ is residually nilpotent, where $N$ is the normal closure of the set $\{[t^{-k}x_it^{k},x_j]: i,j\in\{1,\ldots, n\}; i\neq j; k\in\Z\}$ in $G$.  Let $g\in G/N$ with $g\neq 1$. If the exponent sum of $t$ in $g$ is non-zero, then we can take the homomorphism $\phi : G/N \ar \Z$ with $x_i\mapsto 0$ and $t\mapsto 1$. Then $\phi(g)\neq 1$ and since $\Z$ is residually nilpotent the result follows.
On the other hand, assume that the exponent sum of $t$ in $g$ is zero.  Notice that the relations $[t^{-k}x_it^{k},x_j]$ are equivalent to $[t^{-s}x_it^{s},t^lx_jt^{-l}]$ for every $s,l\in\Z$.  Using these relations,  $g$ can be written as
$$
g=(t^{-s_1}x_1^{w_1}t^{s_1}\cdots t^{-s_{m_1}}x_1^{w_{m_1}}t^{s_{m_1}}) \cdots (t^{-q_1}x_k^{z_1}t^{q_1}\cdots t^{-q_{m_k}}x_k^{z_{m_k}}t^{q_{m_k}}),
$$
where the words 
$$
(t^{-s_1}x_1^{w_1}t^{s_1}\cdots t^{-s_{m_1}}x_1^{w_{m_1}}t^{s_{m_1}}), \ldots,  (t^{-q_1}x_k^{z_1}t^{q_1}\cdots t^{-q_{m_k}}x_k^{z_{m_k}}t^{q_{m_k}})
$$ 
are reduced. Since $g\neq 1$ at least one of the $w_1,\ldots, w_{m_1},\ldots, z_1,\ldots, z_{m_k}\neq 0$. Then we take the homomorphism $\phi : G/N\ar \Z*\Z_{p_i^{r_i}}$ with $\phi(\<t\>)=\Z$, $\f(\<x_i\>)=\Z_{p_i^{r_i}}$ and $\f(x_j)=0$ for all $j\neq i$. Since $\Z*\Z_{p_i^{r_i}}$ are residually nilpotent (see \cite[Theorem 4.1]{gruenberg}),   the result follows.\qed

\begin{lemma}\label{new-zzm} 
For a positive integer $m$, with $m \geq 2$, let $G=\<t,a\mid a^m=1\>$. Then $\gamma_{\w}(G)$ is the normal closure of the set $\{[t^{-k}a^{\mu}t^{k},a^{\nu}]: k \in \mathbb{Z}; \mu, \nu \in \mathbb{N}; {\rm gcd}(\mu, \nu) = 1; \mu \nu = m\}$ in $G$.
\end{lemma}

\pf Let $m=p_1^{r_1}\cdots p_n^{r_n}$, with $n\ge 2$ be the prime number decomposition of $m$. For $i \in \{1, \ldots, n\}$, let $q_{i} = \frac{m}{p^{r_{i}}_{i}}$. Since ${\rm gcd}(q_{1}, \ldots, q_{n}) = 1$, there are $d_{i} \in \mathbb{Z}$ such that $q_{1}d_{1} + \cdots + q_{n}d_{n} = 1$. For $i \in \{1, \ldots, n\}$, let $u_{i} = a^{q_{i}d_{i}}$. Then, for any $i \in \{1, \ldots, n\}$, the order of $u_{i}$ is $p^{r_{i}}_{i}$, and $a = u_{1}u_{2} \cdots u_{n}$.  Since $G\cong \Z*(\Z_{p_1^{r_1}}\times\cdots\times\Z_{p_n^{r_n}})$,  $G$ admits a presentation as in  Lemma \ref{zzm} and the isomorphism between the two presentations implies that each $x_i$ maps to $u_i$. 

Let $\mathcal{M}=\{[t^{-k}u_it^k,u_j] ; i,j\in \{1,\ldots,n\}; i\neq j; k\in\Z\}$ be the generating set of $\gamma_{\omega}(G)$ described in Lemma \ref{zzm} and $\mathcal{K} = \{[t^{-k}a^{\mu}t^{k},a^{\nu}]: k \in \mathbb{Z}; \mu, \nu \in \mathbb{N}; {\rm gcd}(\mu, \nu) = 1; \mu \nu = m\}$. Write $K$ for the normal closure of the set $\mathcal{K}$ in $G$.  We claim that $\gamma_{\omega}(G) = K$.  We first show that $\gamma_{\omega}(G) \subseteq K$. Let $[t^{-k}u_{i}t^{k},u_{j}] \in {\mathcal{M}}$, and, without loss of generality, we assume that $i < j$. Write $s_{1} = p^{r_{i+1}}_{i+1} \cdots p^{r_{j-1}}_{j-1}p^{r_{j+1}}_{j+1} \cdots p^{r_{n}}_{n}$ and $s_{2} = p^{r_{1}}_{1} \cdots p^{r_{i-1}}_{i-1}$.  Then  
$$
[t^{-k}u_{i}t^{k},u_{j}] = [t^{-k}a^{q_{i}d_{i}}t^{k},a^{q_{j}d_{j}}] = [t^{-k}(a^{p_1^{r_1}\ldots p_{i-1}^{r_{i-1}}p_j^{r_j}})^{d_i s_1}t^{k}, (a^{p_i^{r_i}\cdots p_{j-1}^{r_{j-1}}p_{j+1}^{r_{j+1}}\cdots p_n^{r_n}})^{d_j s_2}].
$$ 
Working in $G/K$,  
$$
[t^{-k}(a^{p_1^{r_1}\cdots p_{i-1}^{r_{i-1}}p_j^{r_j}})t^{k}, (a^{p_i^{r_i}\cdots p_{j-1}^{r_{j-1}}p_{j+1}^{r_{j+1}}\cdots p_n^{r_n}})]=1
$$ 
for all $k \in \mathbb{Z}$. By using the commutator identities $[xy,z]=[x,z]^y[y,z]$, $[x,yz]=[x,z][x,y]^z$, $[x^{-1},y] = ([x,y]^{-1})^{x^{-1}}$ and $[x,y^{-1}] = ([x,y]^{-1})^{y^{-1}}$ repeatedly, we get 
$$
[t^{-k}(a^{p_1^{r_1}\cdots p_{i-1}^{r_{i-1}}p_j^{r_j}})^{d_i s_1}t^{k}, (a^{p_i^{r_i}\cdots p_{j-1}^{r_{j-1}}p_{j+1}^{r_{j+1}}\ldots p_n^{r_n}})^{d_j s_2}]=1
$$ 
for all $k \in \mathbb{Z}$. Hence $[t^{-k}u_{i}t^{k},u_{j}] \in K$ for $i < j$.  Applying similar arguments as above, we have, for $i > j$, $[t^{-k}u_{i}t^{k},u_{j}] \in K$. Consequently, $\gamma_{\omega}(G)\subseteq K$.

For the converse, since $a = u_{1}u_{2} \cdots u_{n}$, ${\rm gcd}(\mu,\nu) = 1$ and $\mu \nu = m$, the elements of $\mathcal{K}$ are 
$$
[t^{-k}a^{\mu}t^{k},a^{\nu}] = [t^{-k}(u_1\cdots u_n)^{\mu}t^{k},(u_1\cdots u_n)^{\nu}]=[t^{-k} (u_{i_1}\cdots u_{i_{l}})^{\mu}t^{k}, (u_{j_1}\cdots u_{j_{n-l}})^{\nu}],
$$ 
with 
$\{ u_{i_1},\ldots, u_{i_l}\}\bigcup \{u_{j_1},\ldots, u_{j_{n-l}}\}=\{u_1,\ldots, u_n\}$ and $\{ u_{i_1},\ldots, u_{i_{l}}\}\bigcap \{u_{j_1},\ldots, u_{j_{n-l}}\}=\emptyset.$ Now one can easily show,  by using the commutator identities $[xy,z]=[x,z]^y[y,z]$ and $[x,yz]=[x,z][x,y]^z$ repeatedly,  that the elements of $K$ belong to $\gamma_{\omega}(G)$.  Therefore $K \subseteq \gamma_{\omega}(G)$ and so $\gamma_{\omega}(G) = K$.  \qed

\subsection{Known results on Baumslag--Solitar groups}

Moldavanskii in \cite{mold} has shown the following.

\begin{proposition}[{\cite[Theorem 1]{mold}}]\label{Nw}
Let $G=\BS(m,n)$, with $0 < m \leq |n|$ and $d={\rm gcd}(m,n)$. Then $N_{\w}$ coincides with the normal closure of the set $\{[t^ka^dt^{-k},a]: k \in \mathbb{Z}\}$ in $G$.
\end{proposition}

By Proposition \ref{Nw} and Lemma \ref{bounds}, we get the following  result which we will use  repeatedly in the following.

\begin{cor}\label{user}
Let $G=\BS(m,n)$, with $0 < m \leq |n|$ and $d= {\rm gcd}(m,n)$. Then for all $k, x, y \in \mathbb{Z}$, $[(t^{-k}a^{d}t^k)^{x}, a^y] \in \g_{\w}(G)$ and $[(t^{-k}at^k)^{y}, (a^{d})^{x}] \in \g_{\w}(G)$.
\end{cor}

\pf Since $[t^{-k}a^{y}t^k, a^{dx}] = [a^y, (a^{dx})^{t^{-k}}]^{t^{k}}$ and $\gamma_{\w}(G)$ is normal in $G$, it suffices to prove that $[t^{-k}a^{dx}t^k, a^y] \in \g_{\w}(G)$ for all $k, x, y \in \mathbb{Z}$. By Proposition \ref{Nw} and Lemma \ref{bounds}, we get $[t^{-k}a^{d}t^k, a] \in \g_{\w}(G)$ for all $k \in \mathbb{Z}$.  By using a double induction argument on $x$ and $y$, we obtain the desired result. \qed

\vskip .120 in

Moreover, Moldavanskii in \cite{mold-ru} (see also \cite{mold-rev}) has shown the following.

\begin{proposition}\label{Np}
Let $G=\BS(m,n)$, $p$ be a prime number and let $m=p^rm_1$ and $n=p^sn_1$ where $r,s\ge 0$ and $m_1,n_1$ are not divided by $p$. Let also $d={\rm gcd}(m_1,n_1)$, $m_1=du$ and $n_1=dv$. Then
\begin{enumerate}
\item if $r\neq s$ or if $m_1\not\equiv n_1 \pmod p$, then $(Np)_{w}$ coincides with the normal closure of $a^{p^\xi}$ in $G$,  where $\xi=\min\{r,s\}$.
\item if $r=s$ and $m_1\equiv n_1 \pmod p$, then $(Np)_{w}$ coincides with the normal closure  of the set $\{t^{-1}a^{p^ru}ta^{-p^{r}v}, [t^ka^{p^r}t^{-k},a]: k\in{\mathbb Z}\}$ in $G$.
\end{enumerate}
\end{proposition}

\section{Calculation of $\g_{\w}(\BS(m,n))$}

\begin{proposition}\label{coprime}
Let $G=\BS(m,n)$, with $0 < m \leq |n|$ and ${\rm gcd}(m,n)=1$.  Then 
\begin{enumerate}
\item If there is a prime $p$ such that $n\equiv m \pmod p$,  then $\g_{\w}(G)$ is the normal closure of the set $\{[t^{-k}at^{k},a]: k\in\Z\}$ in $G$.

\item If $n\not\equiv m \pmod p$ for any prime integer $p$, then $\g_{\w}(G)$ is the normal closure of $a$ in $G$.

\end{enumerate}
\end{proposition}

\pf \begin{enumerate}
\item Assume that there is a prime number $p$ such that $n\equiv m \pmod p$. Since ${\rm gcd}(m,n)=1$, we have $p$ divides neither $m$ nor $n$ and so $m$ and $n$ satisfy the conditions of Proposition \ref{Np} (2). Therefore we have $(Np)_{w}$ is the normal closure in $G$ of $[t^{-k}at^k,a]$. On the other hand, by Proposition \ref{Nw}, we have $N_{\w}$ is the normal closure in $G$ of the set $\{[t^{-k}at^{k},a]: k\in\Z\}$. Hence the description of $\g_{\w}(G)$ is an immediate consequence of  Lemma \ref{bounds}. 

\item  Assume that $m\not\equiv n \pmod p$ for every prime number $p$.  Since $\gcd(m,n)=1$, we have by Proposition \ref{Np}~(1),  $(Np)_{\w}$ is the normal closure of $a$ in $G$. So $a \in \bigcap\limits_{p\ {\rm prime}}(Np)_{w}$. By Lemma \ref{bounds},  we obtain the required result.   \qed
\end{enumerate}

\begin{remark}\label{rem1}\upshape{Notice that in the above Proposition,  when $m=1$ and $n\neq 2$, then the commutators $[t^{-k}at^k,a]$ are trivial and hence $\g_{\w}(G)=\{1\}$.
}
\end{remark}

\begin{lemma}\label{kofinas}
Let $G=\BS(m,n)$, with $0 < m \leq |n|$, let ${\rm gcd}(m,n)=d$ and let $\mu,\nu$ be positive integers such that $1\le\mu,\nu\le d$, ${\rm gcd}(\mu,\nu)=1$ and $\mu\nu=d$. Then $[t^{-k}a^{\mu}t^{k},a^{\nu}] \in \gamma_{\w}(G)$ for all $k \in \mathbb{Z}$. 
\end{lemma}

\pf By Corollary \ref{user}, we have $[t^{-k}a^{d}t^{k}, a] \in \g_{\w}(G)$. Hence, in the case $d$ is a power of a prime number, that is, $\mu = 1$ or $\nu = 1$, the required result follows. Thus, in what follows, we may assume that $\mu, \nu > 1$. Fix some $k \in \mathbb{Z}$ and let us denote $u=t^{-k}a^{\mu}t^{k}$. Assume that $[u, a^{\nu}] \in \gamma_{i}(G)$ for some $i \geq 2$. Since $[u^{\nu}, a^{\nu}] = [t^{-k}a^{\mu\nu}t^{k}, a^{\nu}]=[t^{-k}a^dt^{k}, a^{\nu}]$, it follows from Corollary \ref{user} that $[u^{\nu}, a^{\nu}] \in\g_{\w}(G)$ and hence, $[u^{\nu}, a^{\nu}] \in \g_j(G)$ for all $j \in {\mathbb N}$. In particular, we have $[u^{\nu}, a^{\nu}] \in \g_{i+1}(G)$. Since $[u, a^{\nu}] \in \gamma_{i}(G)$ and $[u^{\nu}, a^{\nu}] \in \g_{i+1}(G) = [\gamma_{i}(G),G]$, we get, by Lemma \ref{rep5}~(3) (for $N = \gamma_{i}(G)$), $[u, a^{\nu}]^{\nu} \in \g_{i+1}(G)$. 
Similarly, since $[u, a^{\mu\nu}] = [u, a^{d}] = [a^{\mu}, t^{k}a^{d}t^{-k}]^{t^{k}}$ and $\gamma_{\w}(G)$ is normal in $G$, it follows from Corollary \ref{user} that $[u, a^{\mu\nu}] \in \g_{\w}(G)$. In particular, we have $[u, a^{\mu\nu}] \in\g_{i+1}(G)$. As before, by Lemma \ref{rep5}~(2), we get $[u, a^{\nu}]^{\mu} \in \g_{i+1}(G)$. Thus, $[u, a^{\nu}]^{\nu}, [u, a^{\nu}]^{\mu} \in \g_{i+1}(G)$. Since ${\rm gcd}(\mu, \nu) = 1$, we have $[u, a^{\nu}] \in \g_{i+1}(G)$. We carry on this process  and we obtain the required result.  \qed

We are now able to give the proof of our main theorem.

\setcounter{theorem}{0}
\begin{theorem}
Let $G=\BS(m,n)$, with $0 < m \leq |n|$, ${\rm gcd}(m,n)=d \geq 2$, $m=dm_1$ and $n=dn_1$.   
Then  
\begin{enumerate}
\setcounter{theorem}{2}

\item If $n_1\not\equiv m_1 \pmod p$ for every prime number $p$, then $\g_{\w}(G)$ is the normal closure  of the set $\{a^{d}, [t^{-k}a^{\mu}t^{k},a^{\nu}]: k\in\Z, \mu,\nu \in \mathbb{N}, {\rm gcd}(\mu, \nu) = 1 ~{\rm and}~\mu \nu = d\}$ in $G$.  

\item If there is a prime number $p$ such that  $n_1\equiv m_1\pmod p$, then $\g_{\w}(G)$ is the normal closure of the set $\{[t^{-k}a^{\mu}t^{k},a^{\nu}]: k\in\Z, \mu,\nu \in \mathbb{N}, {\rm gcd}(\mu, \nu) = 1 ~{\rm and}~\mu \nu = d\}$  in $G$. 

\end{enumerate}
\end{theorem}

\begin{remark}\upshape{Notice that if $n_1\not\equiv m_1\pmod p$ for every prime $p$, then $n_1-m_1\not\equiv 0 \pmod p$ for every prime $p$, therefore $n_1-m_1=\pm 1$.  Hence the two possibilities of Theorem \ref{main} can be simplified as to whether $m_1=n_1\pm 1$ or not.
}
\end{remark}

\pf
\begin{enumerate}

\item Let $T$ be the normal closure  of the set $\{a^{d}, [t^{-k}a^{\mu}t^{k},a^{\nu}]: k\in\Z, \mu,\nu \in \mathbb{N}, {\rm gcd}(\mu, \nu) = 1 ~{\rm and}~\mu \nu = d\}$ in $G$.  We first show that $T \subseteq \gamma_{\w}(G)$. By Lemma \ref{kofinas}, it is enough to show that $a^{d} \in \gamma_{\w}(G)$. Let $d = p^{r_{1}}_{1} \cdots p^{r_{\kappa}}_{\kappa}$ be the prime number decomposition of $d$. For $i \in \{1, \ldots, \kappa\}$, we write $m = p^{r_{i}}_{i} m^{\prime}_{1i}$ and $n = p^{r_{i}}_{i} n^{\prime}_{1i}$, where $m^{\prime}_{1i} = \frac{d}{p^{r_{i}}_{i}}m_{1}$ and $n^{\prime}_{1i} = \frac{d}{p^{r_{i}}_{i}}n_{1}$. Now, $n^{\prime}_{1i} - m^{\prime}_{1i} = \frac{d}{p^{r_{i}}_{i}}(n_{1} - m_{1})$. Since $p_{i} \nmid (n_{1}-m_{1})$, we obtain $p_{i} \nmid (n^{\prime}_{1i}-m^{\prime}_{1i})$. By Proposition \ref{Np}~(1), we have $a^{p^{r_{i}}_{i}} \in (Np_{i})_{\w}$ for every $i \in \{1, \ldots, \kappa\}$ and so, $a^{d} \in (Np_{i})_{\w}$ for any $i \in \{1, \ldots, \kappa\}$. By Proposition \ref{Np}~(1), for every $q \notin \{p_{1}, \ldots, p_{\kappa}\}$, we have $a \in (Nq)_{\w}$, which again implies that $a^{d} \in (Nq)_{\w}$. Therefore $a^{d} \in \bigcap\limits_{p\ {\rm prime}}(Np)_{w}$. By Lemma \ref{bounds}, we have $a^{d} \in \gamma_{\w}(G)$ and so $T \subseteq \gamma_{\w}(G)$.  On the other hand,  by Lemma \ref{new-zzm}, $G/T\cong (\Z*\Z_d)/\g_{\w}(\Z*\Z_d)$ which is residually nilpotent. Hence $\g_{\w}(G)\subseteq T$. Therefore $\g_{\w}(G)=T$ and we obtain the required result.

\item Assume that there is some prime $p$ such that $n_1\equiv m_1 \pmod p$ and let $M$ be the normal closure of the set $\{[t^{-k}a^{\mu}t^{k},a^{\nu}]: k\in\Z, \mu,\nu \in \mathbb{N}, {\rm gcd}(\mu, \nu) = 1 ~{\rm and}~\mu \nu = d\}$ in $G$.    By Lemma \ref{kofinas}, we have $M \subseteq \g_{\w}(G)$. We claim that $G/M$ is residually nilpotent. 
Notice that $G/M$ has a presentation of the form
$$G/M=\<a,t\mid t^{-1}a^mt=a^n, [t^{-k}a^{\mu}t^{k},a^{\nu}]=1, \ k\in\Z\>$$
with $ \mu,\nu \in \mathbb{N}, {\rm gcd}(\mu, \nu) = 1 ~{\rm and}~\mu \nu = d.$
Let $g$ be an element in $G/M$. Then $g$ is a word in $a,t$.  

Assume first that the exponent sum of $t$ in $g$ is non-zero. Then there is a map $\f:G/M\ar \Z$ such that $a\mapsto 0$ and $t\mapsto 1$. It can easily be seen that $\f$ is a homomorphism and that $\f(g)\neq 1$.  Since $\Z$ is nilpotent, the result follows.

Assume now that the exponent sum of $t$ in $g$ is zero.  Then $g$ can be written in reduced form
$$
g = a^{\rho_0}(t^{\varepsilon_{1}}a^{\rho_1}t^{-\varepsilon_{1}})(t^{\varepsilon_{1}+\varepsilon_{2}}a^{\rho_2}t^{-(\varepsilon_{1}+\varepsilon_{2})}) \cdots (t^{\varepsilon_{1}+ \cdots +\varepsilon_{\kappa-1}}a^{\rho_{\kappa-1}}t^{-(\varepsilon_{1}+ \cdots + \varepsilon_{\kappa-1})})a^{\rho_{\kappa}} \eqno(3.1)
$$
with $|\rho_i|<m$ if $\e_1+\cdots+\e_i\le-1$ and $|\rho_i|<n$ if $\e_1+\cdots+\e_i\ge1$.  Write each $a^{\rho_i} = a^{d \lambda_i} a^{r_i}$, with $r_i \in \{0, \ldots, d-1\}$, $t^{-k} a^{\rho_i}t^{k} = (t^{-k}(a^{d})^{\lambda_i}t^{k})(t^{-k}a^{r_i}t^{k})$. Using the identity $t^{-k}a^{\zeta}t^{k} = a^{\zeta}[a^{\zeta},t^{k}]$, we rewrite all the above and replace them in the expression $(3.1)$. Then, using the identity $ab=ba[a,b]$ as many times as needed and the identities 
$$
[ab,c] = [a,c]~[[a,c],b]~[b,c], \eqno(3.2)
$$ 
$$
[a,bc] = [a,c]~[a,b]~[[a,b],c], \eqno(3.3)
$$ 
$g$ has an expression of the form
$$
g=a^{\lambda}~[(a^{d})^{\lambda_{1}},t^{k_1}] \ldots [(a^{d})^{\lambda_{s}},t^{k_s}]\cdot w
$$
where $\lambda_{1}, \ldots, \lambda_{s} \in \mathbb{N}$ and $w$ is a product of group commutators of the form $[h_{1}, \ldots, h_{r}]$, with $r \geq 2$ and $h_{1}, \ldots, h_{r} \in \{a, \ldots, a^{d-1}\} \cup \{t^{k}: k \in \mathbb{Z} \setminus \{0\}\}$. Note that $d\lambda_{1}, \ldots, d\lambda_{s} < m_{1}, n_{1}$.  Next, we separate two cases.
\begin{enumerate}

\item Let $w=1$. For the next few lines, let $G_{1} = {\rm BS}(m_{1},n_{1}) = \langle \bar{t}, \bar{a}: (\bar{t})^{-1}(\bar{a})^{m_{1}}\bar{t} = (\bar{a})^{n_{1}} \rangle$. Since ${\rm gcd}(m_{1}, n_{1}) = 1$ and $n_{1} \equiv m_{1} \pmod p$ for a prime integer $p$, we have, by Proposition \ref{coprime}~(1) and $(3.2)$, $\gamma_{\w}(G_{1})$ is the normal closure  of the set $\{[\bar{t}^{-k}\bar{a}\bar{t}^{k}, \bar{a}]: k \in \mathbb{Z}\}$ in $G_{1}$.  Then there is a natural homomorphism $\psi: G/M\ar G_1/\g_{\w}(G_1)$ with $a\mapsto \bar{a}$ and $t\mapsto\bar{t}$ such that $\psi(g)\neq 1$. Since $G_1/\g_{\w}(G_1)$ is residually nilpotent, the result follows.

\item Let $w\neq 1$.  Then it suffices to map $G/M$ to $(G/M)/\<a^d\>\cong (\Z*\Z_d)/\g_{\w}(\Z*\Z_d)$. The image of $g$ is $w\neq 1$ which is reduced.  The result follows from that fact that $(G/H)/\<a^d\>\cong (\Z*\Z_d)/\g_{\w}(\Z*\Z_d)$ is residually nilpotent.  \qed
\end{enumerate} 
\end{enumerate}

\section{The group $[\g_{\w}(G),G]$}

\begin{lemma}\label{lekmp2}
Let $G = {\rm BS}(m, n)$, with $0 < m \leq |n|$, and let $d = {\rm gcd}(m, n)$. Then 
\begin{enumerate}
\item $[t^{-k}a^{d}t^{k}, a^{d}] \in [\g_{\w}(G), G]$ for all $k \in \mathbb{Z}$.

\item $[t^{-k}a^{d}t^{k}, a]^{d} \in [\g_{\w}(G), G]$ for all $k \in \mathbb{Z}$.

\item $[t^{-k}at^{k}, a^{d}]^{d} \in [\g_{\w}(G), G]$ for all $k \in \mathbb{Z}$.
\end{enumerate}
\end{lemma}

\pf Throughout the proof, we write $u_k = t^{-k}a^dt^k$, with $k \in \mathbb{Z}$. By Corollary \ref{user}, $[u_k, a^{d}] \in \g_{\w}(G)$ for all $k \in \mathbb{Z}$.  Furthermore, we write $m = dm_{1}$ and $n = dn_{1}$, where $\gcd(m_{1}, n_{1}) = 1$. 
\begin{enumerate}
\item Since $[u_k ,a^d]=([u_{-k}, a^d]^{-1})^{t^k}$ and $[\g_{\w}(G),G]$ is normal in $G$, it suffices to show that $[u_k, a^{d}] \in [\g_{\w}(G), G]$ for all $k \in \mathbb{N}$. We use induction on $k$. 
Assume at first that $k = 1$. Since, as aforementioned, $[u_1, a^{d}] \in \g_{\w}(G)$ and since $[u_{1}^{m_{1}}, a^{d}] = [t^{-1}a^{m}t, a^{d}] = [a^{n},a^{d}] = 1$, it follows from Lemma \ref{rep5} (3) (for $N=\g_{\w}(G)$) that $[u_1, a^{d}]^{m_{1}} \in [\g_{\w}(G), G]$.  
Since $[u_1, a^{dn_{1}}] = [a^{d}, ta^{n}t^{-1}]^{t} = [a^{d}, a^{m}]^{t} = 1 in [\g_{\w}(G),G]$, by Lemma \ref{rep5} (2)  (for $N=\g_{\w}(G)$) we get $[u_1, a^{d}]^{n_{1}} \in [\g_{\w}(G), G]$.  But $\gcd(m_{1}, n_{1}) = 1$ and so, the result follows for $k=1$.

Assume that  $[u_k, a^{d}] \in [\g_{\w}(G), G]$ for some $k \in \mathbb{N}$.  Then, by Lemma \ref{rep5} (1) (for $N=[\g_{\w}(G),G]$) we have $[u_k^x,a^d]\in [\g_{\w}(G),G]$ for any $x\in\N$. 
Hence, 
$$
[u_{k+1}^{m_{1}}, a^{d}] = [t^{-(k+1)}a^mt^{k+1},a^d] = [t^{-k}(t^{-1}a^{m}t)t^{k}, a^{d}] = [u_k^{n_{1}},a^d] \in [\g_{\w}(G),G]
$$
Since $[u_{k+1},a^d]\in\g_{\w}(G)$, it follows from Lemma \ref{rep5} (2) (for $N=\g_{\w}(G)$) that $[u_{k+1}, a^d]^{m_1}\in [\g_{\w}(G),G]$. 

As above,  since $[u_k,a^d]\in[\g_{\w}(G),G]$, we have,  by Lemma \ref{rep5} (1) (for $N=[\g_{\w}(G),G]$ that $[u_k,a^{dm_1}]\in [\g_{\w}(G),G]$ and therefore 
$$
[u_{k+1}, a^{dn_{1}}]^{t^{-1}} = [t^{-(k+1)}a^{d}t^{k+1},a^n]^{t^{-1}} = [u_{k}, ta^{n}t^{-1}] = [u_{k}, a^{m}] = [u_k, a^{dm_{1}}] \in [\g_{\w}(G),G]
$$
But $[\g_{\w}(G),G]$ is normal in $G$,  so $[u_{k+1}, a^{dn_{1}}] \in [\g_{\w}(G),G]$.  Again, by Lemma \ref{rep5} (2) (for $N=\g_{\w}(G)$),  $[u_{k+1}, a^{d}] ^{n_1}\in [\g_{\w}(G),G]$. Since $\gcd(m_{1}, n_{1}) = 1$ we obtain the required result.

\item By Corollary \ref{user}, $[u_k, a] \in \g_{\w}(G)$ and, by Lemma \ref{lekmp2} (1), $[u_k, a^d] \in [\g_{\w}(G), G]$ for all $k \in \mathbb{Z}$, the result follows from Lemma \ref{rep5} (2) (for $N=\g_{\w}(G)$.

\item By Corollary \ref{user}, $[t^{-k}at^{k}, a^{d}] \in \g_{\w}(G)$. Since $[t^{-k}at^{k}, a^{d}] = ([t^{k}a^{d}t^{-k}, a]^{-1})^{t^{k}}$, the result follows from Lemma \ref{lekmp2} (2). \qed
\end{enumerate}

\begin{proposition}\label{prokmp1}
Let $G = {\rm BS}(m, n)$, with $0 < m \leq |n|$ and let $d = \gcd(m, n)$ be a power of a prime integer $p$. Then $[t^{-k}a^{d}t^{k}, a], [t^{-k}at^{k}, a^{d}] \in [\g_{\w}(G), G]$ for all $k \in \mathbb{N}$.
\end{proposition}

\pf 
Let $d = p^{\mu}$, with $\mu \geq 1$. Thus we may write $m = p^{r}m_{1}$ and $n = p^{s}n_{1}$, where $\mu = \min\{r, s\}$ and $\gcd(p, m_{1}) = \gcd(p, n_{1}) = 1$. Throughout the proof, we write $u_k = t^{-k}a^{d}t^k$ and $v_{k} = t^{-k}at^{k}$. By Corollary \ref{user}, we have $[u_{k}, a], [v_k, a^d] \in \g_{\w}(G)$ for all $k \in \mathbb{N}$. We separate several cases. In the following, we repeatedly use the fact that $[\g_{\w}(G),G]$ is normal in $G$.

\begin{enumerate}
\item Let $r = s$. In this case, we have $d = p^{r}$. At first, we show that $[u_k, a] \in [\g_{\w}(G), G]$ for all $k \in \mathbb{N}$. We use induction on $k$. Let $k = 1$. Since $[u_{1}^{m_{1}}, a] = [t^{-1}a^{m}t, a] = [a^{n},a] = 1 \in [\gamma_{\w}(G),G]$ and $[u_{1},a] \in \gamma_{\w}(G)$,  we have from Lemma \ref{rep5} (3) (for $N=\g_{\w}(G))$ that $[u_1, a]^{m_{1}} \in [\g_{\w}(G), G]$. Furthermore, by Lemma \ref{lekmp2}~(2), $[u_1, a]^{p^{r}} \in [\g_{\w}(G), G]$. 
But, $\gcd(m_{1}, p^{r}) = 1$ and so we have $[u_1, a] \in [\g_{\w}(G), G]$. Thus our claim is valid for $k = 1$. 
Assume that $[u_k, a] \in [\g_{\w}(G), G]$ for some $k \in \mathbb{N}$. Using $(3.2)$ as many times as needed and since $[u_k, a] \in [\g_{\w}(G), G]$, we get $[u_{k}^{n_{1}}, a] \in [\g_{\w}(G), G]$. Since $[u_{k+1}^{m_{1}}, a] = [t^{-(k+1)}a^{m}t^{k+1}, a] = [t^{-k}a^{n}t^{k}, a] = [u_{k}^{n_{1}}, a] \in [\g_{\w}(G), G]$ and $[u_{k+1},a] \in \gamma_{\w}(G)$, it follows from Lemma \ref{rep5} (3) (for $N=\g_{\w}(G))$ that $[u_{k+1}, a]^{m_{1}} \in [\g_{\w}(G), G]$. Furthermore, by Lemma \ref{lekmp2}~(2), $[u_{k+1}, a]^{p^{r}} \in [\g_{\w}(G), G]$. But $\gcd(m_{1}, p^{r}) = 1$ and so we have $[u_{k+1}, a] \in [\g_{\w}(G), G]$. Therefore $[u_{k}, a] \in [\gamma_{\w}(G),G]$ for all $k \in \mathbb{N}$. By $(3.2)$, we get $[u_{k}, a] = [a^{p^{r}}, t^{k}, a] \in [\gamma_{\w}(G),G]$ for all $k \in \mathbb{N}$. 

Next we show that $[v_k, a^{p^{r}}] \in [\g_{\w}(G), G]$ for all $k \in \N$. As before, we use induction on $k$. Let $k = 1$. 
Since $ [v_1, a^{p^{r}n_{1}}] = [v_1, a^{n}] = [a,ta^{n}t^{-1}]^{t} = [a,a^{m}]^{t} =1 \in [\gamma_{\w}(G),G]$ and $[v_{1}, a^{p^{r}}] \in \gamma_{\w}(G)$, it follows from Lemma \ref{rep5} (2) (for $N=\g_{\w}(G))$ that $[v_1, a^{p^{r}}]^{n_{1}} \in [\g_{\w}(G), G]$. 
Furthermore, by Lemma \ref{lekmp2} (3), $[v_1, a^{p^{r}}]^{p^{r}} \in [\g_{\w}(G), G]$. But $\gcd(p^{r}, n_{1}) = 1$ and so we have $[v_1, a^{p^{r}}] \in [\g_{\w}(G), G]$. Thus our claim is true for $k = 1$. Assume that $[v_k, a^{p^{r}}] \in [\g_{\w}(G), G]$ for some $k \in \mathbb{N}$.  
Using $(3.3)$ as many times as needed and since $[v_{k},a] \in [\gamma_{\w}(G),G]$, we get $[v_{k},(a^{p^{r}})^{m_{1}}] = [v_{k}, a^{m}] \in [\g_{\w}(G), G]$. Since $[t^{-(k+1)}at^{k+1}, a^{n}] = [t^{-k}at^{k}, a^{m}]^{t} = [v_{k},a^{m}]^{t}$ and $[\gamma_{\w}(G),G]$ is normal in $G$, we get $[v_{k+1}, a^{n}] = [v_{k+1},(a^{p^{r}})^{n_{1}}] \in [\g_{\w}(G), G]$. Since $[v_{k+1},a^{p^{r}}] \in \gamma_{\w}(G)$, we have from Lemma \ref{rep5} (1) (for $N=\g_{\w}(G))$ that $[v_{k+1}, a^{p^{r}}]^{n_{1}} \in [\g_{\w}(G), G]$. 
Furthermore, by Lemma \ref{lekmp2}~(3), $[v_{k+1}, a^{p^{r}}]^{p^{r}} \in [\g_{\w}(G), G]$. 
But $\gcd(n_{1}, p^{r}) = 1$ and so we have $[v_{k+1}, a^{p^{r}}] \in [\g_{\w}(G), G]$. 
Therefore $[v_{k},a^{p^{r}}] \in [\gamma_{\w}(G),G]$ for all $k \in \mathbb{N}$. By $(3.2)$, we get $[v_{k},a^{p^{r}}] = [a, t^{k}, a^{p^{r}}] \in [\gamma_{\w}(G), G]$ for all $k \in \mathbb{N}$.

\item Let $r < s$.  In this case, we have $d = p^{r}$. By similar arguments as in case (1), we get $[u_k, a] \in [\gamma_{\w}(G), G]$ for all $k \in \mathbb{N}$. Thus it remains to show that $[v_k, a^{p^{r}}] \in [\gamma_{\w}(G), G]$ for all $k \in \mathbb{N}$. By Corollary \ref{user} (for $y = x = 1$ and $d = p^{r}$), $[v_{k}, a^{p^{r}}] \in \gamma_{\w}(G)$ for all $k \in \mathbb{N}$. By Lemma \ref{lekmp2}~(3) (for $d = p^{r}$), $[v_{k}, a^{p^{r}}]^{p^{r}} \in [\gamma_{\w}(G),G]$ for all $k \in \mathbb{N}$. We separate two cases.   

\begin{enumerate}

\item Let $2r \leq s$ and fix a positive integer $k \geq 1$. By Lemma \ref{rep5}~(1) (for $N = \gamma_{\w}(G)$, $x = v_{k+1}$, $y = a^{p^{r}}$, $\kappa = p^{r} p^{s-2r}n_{1}$), we have 
$$
[v_{k+1}, (a^{p^{r}})^{p^{r}p^{s-2r}n_{1}}] \equiv [v_{k+1},a^{p^{r}}]^{p^{r} p^{s-2r}n_{1}}~({\rm mod}~[\gamma_{\w}(G),G]). \eqno(4.1)
$$
Since $[v_{k+1},a^{p^{r}}]^{p^{r}} \in [\gamma_{\w}(G),G]$, we get $[v_{k+1},a^{p^{r}}]^{p^{r} p^{s-2r}n_{1}} \in [\gamma_{\w}(G),G]$ and so by $(4.1)$, we have $[v_{k+1},a^{p^{s}n_{1}}] \in [\gamma_{\w}(G),G]$. But $[v_{k+1},a^{p^{s}n_{1}}] = [v_{k},ta^{p^{s}n_{1}}t^{-1}]^{t} = [v_{k},a^{p^{r}m_{1}}]^{t}$. Since $[\gamma_{\w}(G),G]$ is normal in $G$, we get 
$$
[v_{k},a^{p^{r}m_{1}}] \in [\gamma_{\w}(G),G]. \eqno(4.2)
$$
Since $[v_{k}, a^{p^{r}}] \in \gamma_{\w}(G)$, it follows from Lemma \ref{rep5}~(1) (for $N = \gamma_{\w}(G)$, $x = v_{k}$, $y = a^{p^{r}}$, $\kappa = m_{1}$) that 
$$
[v_{k}, a^{p^{r}m_{1}}] \equiv [v_{k},a^{p^{r}}]^{m_{1}} ~({\rm mod}~[\gamma_{\w}(G),G]).
$$
By $(4.2)$, we obtain $[v_{k},a^{p^{r}}]^{m_{1}} \in [\gamma_{\w}(G),G]$. Since $[v_{k},a^{p^{r}}]^{p^{r}} \in [\gamma_{\w}(G),G]$ and ${\rm gcd}(m_{1}, p^{r}) = 1$, we have $[v_{k},a^{p^{r}}] = [t^{-k}at^{k},a^{d}] \in [\gamma_{\w}(G),G]$.

\item Let $2r > s$ and fix a positive integer $k \geq 1$. Since $[v_{k},a^{p^{r}}] \in \gamma_{\w}(G)$, it follows from Lemma \ref{rep5}~(1) (for $x = v_{k}$, $y = a^{p^{r}}$, $\kappa = m_{1}$) that $[v_{k}, a^{p^{r}m_{1}}] \in \gamma_{\w}(G)$ and 
$$
[v_{k},a^{p^{r}m_{1}}] \equiv [v_{k}, a^{p^{r}}]^{m_{1}}~({\rm mod}~[\gamma_{\w}(G),G]).
$$
Since $[\gamma_{\w}(G),G]$ is normal in $G$,
$$
[v_{k}, a^{p^{r}m_{1}}]^{p^{2r-s}} \equiv [v_{k},a^{p^{r}}]^{p^{2r-s}m_{1}}~({\rm mod}~[\gamma_{\w}(G),G]). \eqno(4.3)
$$
Since $[v_{k+1}, a^{p^{s}n_{1}}] = [v_{k}, ta^{p^{s}n_{1}}t^{-1}]^{t} = [v_{k}, a^{p^{r}m_{1}}]^{t}$ and $[v_{k}, a^{p^{r}m_{1}}] \in \gamma_{\w}(G)$, we have 
$$
[v_{k+1}, a^{p^{s}n_{1}}] \equiv [v_{k}, a^{p^{r}m_{1}}]~({\rm mod}~[\gamma_{\w}(G),G]).
$$
By $(4.3)$, we get 
$$
[v_{k+1}, a^{p^{s}n_{1}}]^{p^{2r-s}} \equiv [v_{k}, a^{p^{r}}]^{p^{2r-s}m_{1}}~({\rm mod}~[\gamma_{\w}(G),G]). \eqno(4.4)
$$
Since $[v_{k+1}, a^{p^{r}}] \in \gamma_{\w}(G)$ and $r < s$, it follows from Lemma \ref{rep5}~(1) (for $x = v_{k+1}$, $y = a^{p^{r}}$, $\kappa = p^{s-r}n_{1}$) that 
$$
[v_{k+1}, a^{p^{s}n_{1}}] \equiv [v_{k+1}, a^{p^{r}}]^{p^{s-r}n_{1}}~({\rm mod}~[\gamma_{\w}(G),G]). \eqno(4.5)
$$
By $(4.4)$ and $(4.5)$, we have 
$$
[v_{k}, a^{p^{r}}]^{p^{2r-s}m_{1}} \equiv [v_{k+1}, a^{p^{r}}]^{p^{r}}~({\rm mod}~[\gamma_{\w}(G),G]). \eqno(4.6)
$$
Since $[v_{k+1},a^{p^{r}}]^{p^{r}} \in [\gamma_{\w}(G),G]$, we obtain by $(4.6)$ that $[v_{k}, a^{p^{r}}]^{p^{2r-s}m_{1}} \in [\gamma_{\w}(G),G]$. Since ${\rm gcd}(p^{r}, p^{2r-s}m_{1}) = p^{2r-s}$, we get $[v_{k}, a^{p^{r}}]^{p^{2r-s}} \in [\gamma_{\w}(G),G]$. If $3r \leq 2s$, then, by applying similar arguments as in case $2r \leq s$, we have $[v_{k}, a^{p^{r}}] \in [\gamma_{\w}(G),G]$. If $3r > 2s$, then, by applying similar arguments as in case $2r > s$, we get $[v_{k},a^{p^{r}}]^{p^{3r-2s}} \in [\gamma_{\w}(G),G]$. Since $2r-s > 3r-2s > \cdots$ and since there is $y$ such that $(y+1)r \leq ys$ (for $\frac{r}{s-r} \in \mathbb{N}$, let $y = \frac{r}{s-r}$ and, for $\frac{r}{s-r} \notin \mathbb{N}$, let $y$ be the integral part of $\frac{r}{s-r}$), by continuing this process, we obtain $[v_{k}, a^{p^{r}}] \in [\gamma_{\w}(G),G]$.
\end{enumerate}

\item Let $s < r$. By applying similar arguments as in case $(2)$, we obtain the required result.
\end{enumerate}
By cases (1), (2) and (3), we get the desired result.  \qed

\begin{proposition}\label{thekmp1}
Let $G = {\rm BS}(m, n)$, with $0 < m \leq |n|$ and let $N_{\w}$ be the intersection of all finite index subgroups of $G$. Then $N_{\w} \leq [\gamma_{\w}(G), G]$.
\end{proposition}

\pf 
Let $d = \gcd(m, n)$. Since, by Lemma \ref{bounds}, $N_{\w} \leq \gamma_{\w}(G)$ and, by Proposition \ref{Nw}, $N_{\w}$ coincides with the normal closure  of the set $\{[t^{-k}a^{d}t^{k}, a]: k \in \mathbb{Z}\}$ in $G$, it suffices to show that $[t^{-k}a^{d}t^{k}, a] \in [\gamma_{\w}(G), G]$ for all $k \in \mathbb{Z}$. 
Furthermore, since $[t^{-k}a^xt^k,a^y]=([t^{k}a^yt^{-k},a^x]^{-1})^{t^k}$ and $[\gamma_{\w}(G),G]$ is normal in $G$,  it is enough to show that $[t^{-k}a^{d}t^{k}, a], [t^{-k}at^{k}, a^{d}] \in [\gamma_{\w}(G), G]$ for all $k \in \mathbb{N}$. 
For $d = 1$, the required result follows from Lemma \ref{lekmp2}~(1) and so, from now on, we may assume that $d > 1$. Let $d = p_{1}^{\mu_{1}}\ldots p_{\lambda}^{\mu_{\lambda}}$ be the prime factor decomposition of $d$. 
To prove the result, we use induction on $\lambda$. For $\lambda = 1$, the required result follows from Proposition \ref{prokmp1}. 
Assume that the result is true for some $\lambda \geq 1$ and let $d = p_{1}^{\mu_{1}}\ldots p_{\lambda+1}^{\mu_{\lambda+1}}$. 
Thus we may write $m = p_{1}^{r_{1}}\ldots p_{\lambda+1}^{r_{\lambda+1}}\mu$ and $n = p_{1}^{s_{1}}\ldots p_{\lambda+1}^{s_{\lambda+1}}\nu$, where $\mu_{i} = \min\{r_{i}, s_{i}\}$ and $\gcd(p_{i}, \mu) = \gcd(p_{i}, \nu) = 1$, with $i = 1, \ldots, \lambda+1$. 
For $j \in \{1, \ldots, \lambda+1\}$, let $u = a^{p_{j}^{\mu_{j}}}$ and let $K_{j}$ be the subgroup of $G$ generated by the set $\{u, t\}$.  
For convenience, we write $m_{j} = m/({p_{j}^{\mu_{j}}})$, $n_{j} = n/({p_{j}^{\mu_{j}}})$ and $d_{j} = d/({p_{j}^{\mu_{j}}})$.  
We point out that $K_{j} = {\rm BS}(m_{j}, n_{j})$ (see \cite[Lemma 7.10]{levitt}).      
Since $\gcd(m_{j}, n_{j}) = p_{1}^{\mu_{1}}\ldots p_{j-1}^{\mu_{j-1}}p_{j+1}^{\mu_{j+1}} \ldots p_{\lambda+1}^{\mu_{\lambda+1}}$, by our inductive argument,  we get $[t^{-k}u^{d_{j}}t^{k}, u] \in [\gamma_{\w}(K_{j}), K_{j}] \subseteq [\gamma_{\w}(G), G]$, that is, $[t^{-k}a^{p_{j}^{\mu_{j}}d_{j}}t^{k}, a^{p_{j}^{\mu_{j}}}] \in [\gamma_{\w}(G), G]$ and so $[t^{-k}a^{d}t^{k}, a^{p_{j}^{\mu_{j}}}] \in [\gamma_{\w}(G), G]$. Since $[t^{-k}a^{d}t^{k}, a] \in \gamma_{\w}(G)$, it follows that $[t^{-k}a^{d}t^{k}, a]^{p_{j}^{\mu_{j}}} \in [\gamma_{\w}(G), G]$. 
Therefore, for $j_{1}, j_{2} \in \{1, \ldots, \lambda+1\}$, with $j_{1} \neq j_{2}$, we have $[t^{-k}a^{d}t^{k}, a]^{p_{j_{1}}^{\mu_{j_{1}}}}, [t^{-k}a^{d}t^{k}, a]^{p_{j_{2}}^{\mu_{j_{2}}}} \in [\gamma_{\w}(G), G]$. Since ${\rm gcd}(p_{j_{1}}^{\mu_{j_{1}}}, p_{j_{2}}^{\mu_{j_{2}}}) = 1$, we get $[t^{-k}a^{d}t^{k}, a] \in [\gamma_{\w}(G), G]$ and the result follows. \qed

\begin{cor}
Let $G=\BS(m,n)$, with $0 < m \leq |n|$ and ${\rm gcd}(m,n)=1$.  Then $[\g_{\w}(G),G]=\g_{\w}(G)$.
\end{cor}

\pf Let us assume first that there is some prime number $p$ such that $m\equiv n\pmod p$.  Then by Proposition \ref{coprime}~(1),  $\g_{\w}(G)=N_{\w}$ and so, by Proposition \ref{thekmp1}, $\g_{\w}(G)=[\g_{\w}(G),G]$. On the other hand,  if $m\not\equiv n\pmod p$ for every prime integer $p$, then  $m-n=\pm 1$ and so, the relation $t^{-1}a^mt=a^n$ becomes $t^{-1}a^{n\pm 1}t=a^n$ or equivalently, $t^{-1}a^{n\pm 1}ta^{-(n\pm 1)}=a^{\mp 1}$ or $[t,a^{-(n\pm 1)}]=a^{\mp 1}$.   By Proposition \ref{coprime}~(2), we have $a\in \g_{\w}(G)$. Since $[t,a^{-(n\pm 1)}]=a^{\mp 1}$, we obtain $a\in[\g_{\w}(G),G]$. By Proposition \ref{coprime}~(2), $\gamma_{\w}(G) \subseteq [\gamma_{\w}(G),G]$ and the result follows.  \qed

\begin{proposition}\label{prokmp1-1}
Let $G = {\rm BS}(m, n)$, with $0 < m \leq |n|$, let $d = \gcd(m, n)$ and let $\mu$, $\nu$ be positive integers, with $1 \leq \mu, \nu \leq d$, such that $\gcd(\mu, \nu) = 1$ and $d = \mu\nu$. Then  $[t^{-k}a^{\mu}t^{k}, a^{\nu}] \in [\gamma_{\w}(G), G]$ for all $k \in \mathbb{Z}$. 
\end{proposition}

\pf Fix some $k \in \mathbb{Z}$. By Proposition \ref{Nw}, $[t^{-k}a^{d}t^{k}, a] \in N_{\w}$ and so, by Proposition \ref{thekmp1}, we have $[t^{-k}a^{d}t^{k}, a] = [t^{-k}a^{\mu \nu}t^{k},a] \in [\gamma_{\w}(G), G]$. Using $(3.3)$ as many times as needed and since $[\gamma_{\w}(G),G]$ is a normal subgroup of $G$, we obtain $[t^{-k}a^{\mu\nu}t^{k}, a^{\nu}] \in [\gamma_{\w}(G), G]$. 
By Lemma \ref{kofinas}, $[t^{-k}a^{\mu}t^{k}, a^{\nu}] \in \gamma_{\w}(G)$ and so, by Lemma \ref{rep5} (3), (for $N=\g_{\w}(G))$,  $[t^{-k}a^{\mu}t^{k}, a^{\nu}]^{\nu} \in [\gamma_{\w}(G), G]$. By Proposition \ref{Nw}, Proposition \ref{thekmp1} and since $N_{\w}$ is a normal subgroup of $G$, we get $[t^{-k}at^{k}, a^{d}] \in [\gamma_{\w}(G), G]$. Using $(3.2)$ as many times as  needed and since $[\gamma_{\w}(G),G]$ is a normal subgroup of $G$, we obtain $[(t^{-k}at^{k})^{\mu}, a^{\mu\nu}] = [t^{-k}a^{\mu}t^{k}, a^{\mu \nu}] \in [\gamma_{\w}(G), G]$. 
By Lemma \ref{kofinas}, $[t^{-k}a^{\mu}t^{k}, a^{\nu}] \in \gamma_{\w}(G)$ and so, by Lemma \ref{rep5}, (2) (for $N=\g_{\w}(G))$, we have $[t^{-k}a^{\mu}t^{k}, a^{\nu}]^{\mu} \in [\gamma_{\w}(G), G]$. But $\gcd(\mu, \nu) = 1$, therefore we get the result. \qed

\medskip

So now we may give an answer to the question of Bardakov and Neschadim.

\setcounter{theorem}{1}
\begin{theorem}
Let $G={\rm BS}(m,n)$, with $0 < m \leq |n|$. Then $[\g_{\w}(G),G]=\g_{\w}(G)$.     
\end{theorem}

\pf  Let $d={\rm gcd}(m,n)$, with $m=dm_1$, $n=dn_1$ and ${\rm gcd}(m_1,n_1)=1$.  Let $n = m$. By Theorem \ref{main}~(1), Proposition \ref{prokmp1-1} and $[\gamma_{\w}(G),G] \leq \gamma_{\w}(G)$, we have $[\gamma_{\w}(G),G] = \gamma_{\w}(G)$.  
If there is a prime number $p$ such that $m_1\equiv n_1\pmod p$, then, by Theorem \ref{main}~(2), Proposition \ref{prokmp1-1} and $[\gamma_{\w}(G),G] \leq \gamma_{\w}(G)$, we have the required result. Finally, let $m_1\not\equiv n_1 \pmod p$ for every prime number $p$. By Proposition \ref{prokmp1-1} and Theorem \ref{main}~(2), it is sufficient to show that $a^d\in [\g_{\w}(G),G]$. Since $m_1\not\equiv n_1 \pmod p$ for every prime number $p$, we have $m_1-n_1=\pm 1$.  Hence,  the relation $t^{-1}a^mt=a^n$ implies that $t^{-1}a^{d(n_1 \pm 1)}t=a^{dn_1}$ or equivalently $t^{-1}a^{d(n_1\pm 1)}t=a^{d(n_1\pm 1)}a^{\mp d}$ or $[a^{d(n_1\pm 1)}, t]=a^{\mp d}$.  By Theorem \ref{main}~(2), $a^d\in\g_{\w}(G)$ and so, we obtain the required result.   \qed

\section{Quotient groups}\label{sec5}

Our purpose in this section is to show by means of a Lie algebra method that each quotient group $\gamma_{c}(G)/\gamma_{c+1}(G)$ is finite for any Baumslag-Solitar group $G$.  
Throughout this section, for any group $G$, we write ${\rm gr}_{c}(G) = \gamma_{c}(G)/\gamma_{c+1}(G)$, with $c \in \mathbb{N}$.  
Also, by a Lie algebra we mean a Lie algebra over $\mathbb{Z}$.  
Let ${\rm gr}(G)$  denote  the (restricted) direct sum of the abelian groups ${\rm gr}_{c}(G)$. It is well known that ${\rm gr}(G)$ has the structure of a Lie algebra by defining a Lie multiplication $[a\gamma_{r+1}(G),b\gamma_{s+1}(G)]$ $ = [a,b]\gamma_{r+s+1}(G)$, where $a \gamma_{r+1}(G)$ and $b\gamma_{s+1}(G)$ are the images of the elements $a \in \gamma_{r}(G)$ and $b \in \gamma_{s}(G)$ in the quotient groups ${\rm gr}_{r}(G)$ and ${\rm gr}_{s}(G)$, respectively, and $[a,b]\gamma_{r+s+1}(G)$ is the image of the group commutator $[a,b]$ in the quotient group ${\rm gr}_{r+s}(G)$. Multiplication is then extended to ${\rm gr}(G)$ by linearity.

Let $F$ be a free group of finite rank $n \geq 2$, with a free generating set $\{x_{1}, \ldots, x_{n}\}$. It is well  known that ${\rm gr}(F)$ is a free Lie algebra of rank $n$; freely generated by the set $\{y_{1}, \ldots, y_{n}\}$, where $y_{i} = x_{i}\gamma_{2}(F)$ and $i \in \{1, \ldots, n\}$. Let $N$ be a normal subgroup of $F$. For $c \in \mathbb{N}$, let $N_{c} = N \cap \gamma_{c}(F)$. Note that $N_{1} = N$. Furthermore, for all $c \in \mathbb{N}$, we write ${\rm I}_{c}(N) = N_{c}\gamma_{c+1}(F)/\gamma_{c+1}(F)$. Form the (restricted) direct sum ${\mathcal{L}}(N)$ of the abelian groups ${\rm I}_{c}(N)$. Since $N$ is normal in $F$, ${\mathcal{L}}(N)$ is an ideal of ${\rm gr}(F)$ (see \cite{laz}).

The following result is probably known, but we give a proof for completeness.

\begin{lemma}\label{l5.1}
Let $F$ be a free group of finite rank $n$, with $n \geq 2$, and $N$ be a normal subgroup of $F$. Then, for all $c \in \mathbb{N}$,  ${\rm gr}_{c}(F/N) \cong {\rm gr}_{c}(F)/{\rm I}_{c}(N)$.
\end{lemma}

\pf For all $c \in \mathbb{N}$, $\gamma_{c}(F/N) = \gamma_{c}(F)N/N$. We have the following natural isomorphisms as abelian groups 
$$
\begin{array}{rll}
{\rm gr}_{c}(F/N) & \cong & \frac{\gamma_{c}(F)N}{\gamma_{c+1}(F)N} = \frac{\gamma_{c}(F)(\gamma_{c+1}(F)N)}{\gamma_{c+1}(F)N} \\
& \cong & \frac{\gamma_{c}(F)}{\gamma_{c}(F) \cap (\gamma_{c+1}(F)N)}.
\end{array}
$$
Since $\gamma_{c+1}(F) \subseteq \gamma_{c}(F)$, by the modular law, we have $\gamma_{c}(F) \cap (\gamma_{c+1}(F) N) = \gamma_{c+1}(F)(\gamma_{c}(F) \cap N)$. 
Therefore, for all $c \in \mathbb{N}$, 
$$
{\rm gr}_{c}(F/N) \cong \frac{\gamma_{c}(F)}{\gamma_{c+1}(F)(\gamma_{c}(F) \cap N)} \cong \frac{{\rm gr}_{c}(F)}{\frac{\gamma_{c+1}(F)N_{c}}{\gamma_{c+1}(F)}} = \frac{{\rm gr}_{c}(F)}{{\rm I}_{c}(N)}  
$$ 
as abelian groups in a natural way. \qed  

\vskip .120 in

For the rest of this section, let $F$ be a free group of rank $2$, with a free generating set $\{x, y\}$.

\begin{proposition}\label{pr2}
For a non-zero integer $\kappa$, let ${\mathcal{R}}_{\kappa} = \{x\gamma_{2}(F),$ $\kappa(y\gamma_{2}(F))\}$ and let ${\mathcal{L}}_{\mathcal{R}_{\kappa}}$ be the Lie subalgebra of ${\rm gr}(F)$ generated by the set $\mathcal{R}_{\kappa}$. Then ${\mathcal{L}}_{\mathcal{R}_{\kappa}}$ is free on ${\mathcal{R}}_{\kappa}$ and, for any $c \in \mathbb{N}$, ${\mathcal{L}}_{\mathcal{R}_{\kappa}} \cap {\rm gr}_{c}(F)$ has finite index in ${\rm gr}_{c}(F)$. Moreover, for any $c \in \mathbb{N}$, ${\mathcal{L}}_{\mathcal{R}_{\kappa}} \cap {\rm gr}_{c}(F)$ is spanned by all Lie commutators      of degree $c$ containing $\kappa \overline{y}$. 
\end{proposition} 

\pf Let ${\rm gr}_{1,R}(F)$ be the subgroup of ${\rm gr}_{1}(F)$ generated by the set ${\mathcal{R}}_{\kappa} = \{x\gamma_{2}(F),$ $\kappa(y\gamma_{2}(F))\}$. It is a free abelian group of rank $2$ and has a finite index $|\kappa|$ in ${\rm gr}_{1}(F)$. For $c \in \mathbb{N}$, let ${\mathcal{L}}_{c, \mathcal{R}_{\kappa}} = {\mathcal{L}}_{\mathcal{R}_{\kappa}} \cap {\rm gr}_{c}(F)$. Then ${\mathcal{L}}_{\mathcal{R}_{\kappa}} = \bigoplus_{c \geq 1}{\mathcal{L}}_{c, \mathcal{R}_{\kappa}}$ and each ${\mathcal{L}}_{c, \mathcal{R}_{\kappa}}$ is a free abelian subgroup of ${\rm gr}_{c}(F)$. Observe that 
$$
{\mathbb{Q}} \otimes_{\mathbb{Z}} {\rm gr}(F) = {\mathbb{Q}} \otimes_{\mathbb{Z}} {\mathcal{L}}_{\mathcal{R}_{\kappa}}. \eqno(5.1)
$$
Since ${\rm gr}(F)$ is a free Lie algebra, we get ${\mathbb{Q}} \otimes_{\mathbb{Z}} {\rm gr}(F)$ is a free Lie algebra over $\mathbb{Q}$ on any $\mathbb{Q}$-basis of $\mathbb{Q} \otimes_{\mathbb{Z}} {\rm gr}_{1}(F)$. By $(5.1)$ and since ${\rm gr}(F) = \bigoplus_{c \geq 1}{\rm gr}_{c}(F)$ and ${\mathcal{L}}_{\mathcal{R}_{\kappa}} = \bigoplus_{c \geq 1}{\mathcal{L}}_{c, \mathcal{R}_{\kappa}}$, we have, for all $c \in \mathbb{N}$, 
$$
{\mathbb{Q}} \otimes_{\mathbb{Z}} {\rm gr}_{c}(F) = {\mathbb{Q}} \otimes_{\mathbb{Z}} {\mathcal{L}}_{c, \mathcal{R}_{\kappa}}. \eqno(5.2)
$$ 
Since ${\mathcal{L}}_{c, \mathcal{R}_{\kappa}} \leq {\rm gr}_{c}(F)$ for all $c \in \mathbb{N}$, we have by $(5.2)$, for all $c \in \mathbb{N}$, ${\mathcal{L}}_{c, \mathcal{R}_{\kappa}}$ has finite index in ${\rm gr}_{c}(F)$. Write $\overline{x} = x\gamma_{2}(F)$ and $\overline{y} = y\gamma_{2}(F)$. Since ${\rm gr}(F)$ is a free Lie algebra on the set $\{\overline{x}, \overline{y}\}$ and ${\mathcal{L}}_{\mathcal{R}_{\kappa}}$ is generated as a Lie algebra by the set $\{\overline{x}, \kappa \overline{y}\}$, the natural Lie homomorphism $\psi$ from ${\rm gr}(F)$ into ${\mathcal{L}}_{\mathcal{R}_{\kappa}}$, with $\psi(\overline{x}) = \overline{x}$ and $\psi(\overline{y}) = \kappa \overline{y}$, is onto. It is easily verified that, for all $c \in \mathbb{N}$, $\psi$ induces a $\mathbb{Z}$-linear mapping $\psi_{c}$ from ${\rm gr}_{c}(F)$ onto ${\mathcal{L}}_{c,{\mathcal{R}}_{\kappa}}$. Namely, $\psi_{c}(u(\overline{x}, \overline{y})) = u(\overline{x}, \kappa \overline{y})$ for all $u(\overline{x}, \overline{y}) \in {\rm gr}_{c}(F)$. Since ${\mathcal{L}}_{c,{\mathcal{R}}_{\kappa}}$ has finite index in ${\rm gr}_{c}(F)$ and ${\rm gr}_{c}(F)$ is a free abelian group of finite rank, we obtain $\psi_{c}$ is an isomorphism of abelian groups.

Let $\overline{w} \in {\rm Ker}\psi$. Since ${\rm gr}(F)$ is graded, without loss of generality, we may assume that $\overline{w} = w\gamma_{c+1}(F) \in {\rm gr}_{c}(F)$ for some $w \in \gamma_{c}(F)$ and $c \in \mathbb{N}$. To get a contradiction, we assume that $w = w(x,y) \in \gamma_{c}(F) \setminus \gamma_{c+1}(F)$. Then $\psi(\overline{w}) = \psi(w(\overline{x}, \overline{y}) = \psi_{c}(w(\overline{x}, \overline{y})) = w(\overline{x}, \kappa \overline{y})$. By bilinearity of the Lie bracket, we get $\psi_{c}(w(\overline{x}, \overline{y})) = m(\kappa)w(\overline{x}, \overline{y})$, where $m(\kappa) \in \mathbb{Z} \setminus \{0\}$ depends  on $\kappa$. Since $\overline{w} \in {\rm Ker}\psi$, we have $\psi(\overline{w})=m(\kappa) \overline{w} = 0$ in ${\rm gr}_{c}(F)$. That is, $w(x,y)^{m(\kappa)} \in \gamma_{c+1}(F)$, which is a contradiction, since ${\rm gr}_{c}(F)$ is a free abelian group (of finite rank). Therefore, $\psi$ is an isomorphism of Lie algebras. 
Hence ${\mathcal{L}}_{\mathcal{R}_{\kappa}}$ is a free Lie algebra with a free generating set $\mathcal{R}_{\kappa}$. By the elimination theorem (see, for example, \cite[Chapter 2, Section 2.9, Proposition 10]{bour}), ${\mathcal{L}}_{\mathcal{R}_{\kappa}} = \langle \overline{x} \rangle \oplus L(\{\kappa \overline{y}\} \wr \{\overline{x}\})$, where the free Lie algebra $L(\{\kappa \overline{y}\} \wr \{\overline{x}\})$ is the ideal in ${\mathcal{L}}_{\mathcal{R}_{\kappa}}$ generated by $\kappa \overline{y}$. Clearly 
$$
L(\{\kappa \overline{y}\} \wr \{\overline{x}\}) = \langle \kappa \overline{y} \rangle \oplus \left(\bigoplus_{c \geq 2} {\mathcal{L}}_{c,\mathcal{R}_{\kappa}}\right). \eqno(5.3)
$$ 
and so we obtain ${\mathcal{L}}_{c,\mathcal{R}_{\kappa}}$ is spanned by all Lie commutators      of degree $c$ containing $\kappa \overline{y}$. \qed 

\begin{remark}\upshape
The fact that $\psi_c$ is an isomorphism in the above proof can also be shown as follows.
Since ${\rm gr}_{c}(F)$ is a free abelian group of finite rank, we have ${\mathbb{Q}} \otimes_{\mathbb{Z}} {\rm gr}_{c}(F)$ is a finite-dimensional vector space over $\mathbb{Q}$, and any $\mathbb{Z}$-basis of ${\rm gr}_{c}(F)$ may be regarded as a $\mathbb{Q}$-basis of ${\mathbb{Q}} \otimes_{\mathbb{Z}} {\rm gr}_{c}(F)$. Thus $\psi_{c}$ may be extended to a $\mathbb{Q}$-linear mapping $\overline{\psi}_{c}$ from ${\mathbb{Q}} \otimes_{\mathbb{Z}} {\rm gr}_{c}(F)$ onto ${\mathbb{Q}} \otimes_{\mathbb{Z}} {\mathcal{L}}_{c, \mathcal{R}_{\kappa}}$. By $(5.2)$, we obtain $\overline{\psi}_{c}$ is an isomorphism of vector spaces and so, $\psi_{c}$ is an isomorphism of abelian groups.       
\end{remark}

\begin{proposition}\label{pr3}
For $m, n \in \mathbb{Z} \setminus \{0\}$, let $N$ be the normal closure of the element $r = y^{n-m}[x,y^{m}]$ in $F$. Then, for any $c \in \mathbb{N}$, with $c \geq 2$, ${\rm gr}_{c}(F/N)$ is a finite abelian group.
\end{proposition}

\pf Write $\delta = n-m$, $\overline{x} = x\gamma_{2}(F)$ and $\overline{y} = y \gamma_{2}(F)$. Assume that $\delta \neq 0$ and let ${\mathcal{R}}_{\delta} = \{\overline{x}, \delta \overline{y}\}$. By Proposition \ref{pr2}, ${\mathcal{L}}_{\mathcal{R}_{\delta}}$ is a free Lie algebra on the set ${\mathcal{R}}_{\delta}$. By $(5.3)$ (for $\kappa = \delta \neq 0$), ${\rm I}_{1}(N) = \langle \delta \overline{y} \rangle$ and since ${\mathcal{L}}(N)$ is an ideal in ${\rm gr}(F)$, we have $L(\{\delta \overline{y}\} \wr \{\overline{x}\}) \subseteq {\mathcal{L}}(N)$. Hence, for all $c \geq 2$, ${\mathcal{L}}_{c,{\mathcal{R}}_{\delta}} \subseteq {\rm I}_{c}(N)$. By Proposition \ref{pr2}, for any $c \geq 2$, ${\mathcal{L}}_{c,{\mathcal{R}}_{\delta}}$ has finite index in ${\rm gr}_{c}(F)$ and  so ${\rm I}_{c}(N)$ has finite index in ${\rm gr}_{c}(F)$. By Lemma  \ref{l5.1}, we obtain, for any $c \geq 2$, ${\rm gr}_{c}(F/N)$ is a finite abelian group. 

Thus we may assume that $\delta = 0$. In this case let ${\mathcal{R}}_{m} = \{\overline{x}, m \overline{y}\}$. Then ${\rm I}_{1}(N) = \{0\}$. We use similar arguments as before. By Proposition \ref{pr2}, ${\mathcal{L}}_{\mathcal{R}_{m}}$ is a free Lie algebra on the set ${\mathcal{R}}_{\delta}$. By $(5.3)$ (for $\kappa = m$), ${\rm I}_{2}(N) = \langle [\overline{x}, m \overline{y}] \rangle$ and since ${\mathcal{L}}(N)$ is an ideal in ${\rm gr}(F)$, we have $\bigoplus_{c \geq 2}\mathcal{L}_{c,{\mathcal{R}}_{m}} \subseteq {\mathcal{L}}(N)$. Hence, for all $c \geq 2$, ${\mathcal{L}}_{c,{\mathcal{R}}_{m}} \subseteq {\rm I}_{c}(N)$. By Proposition \ref{pr2}, for any $c \geq 2$, ${\mathcal{L}}_{c,{\mathcal{R}}_{m}}$ has finite index in ${\rm gr}_{c}(F)$ and  so ${\rm I}_{c}(N)$ has finite index in ${\rm gr}_{c}(F)$. By Lemma \ref{l5.1}, we obtain, for any $c \geq 2$, ${\rm gr}_{c}(F/N)$ is a finite abelian group.  \qed

\begin{cor}
For $m, n \in \mathbb{Z} \setminus \{0\}$, let $G = {\rm BS}(m,n) = \langle t, a \mid a^{n-m} [t,a^{m}] = 1\rangle$. Then, for any $c \in \mathbb{N}$, with $c \geq 2$, ${\rm gr}_{c}(G)$ is a finite abelian group.
\end{cor}

\pf For $m, n \in \mathbb{Z} \setminus \{0\}$, let $N$ be the normal closure of the element $r = y^{n-m}[x,y^{m}]$ in $F$. Write $t = xN$ and $a = yN$. Clearly the quotient group $F/N$ has a presentation $\langle t, a \mid a^{n-m} [t,a^{m}] = 1\rangle$. So, by Proposition \ref{pr3}, we obtain the desired result. \qed

\end{document}